\documentclass[a4paper]{amsart}
\usepackage{graphicx}
\theoremstyle{plain}
  \newtheorem{thm}{Theorem}[section]
  
  \newtheorem{cor}[thm]{Corollary}

  \newtheorem{claim}[thm]{Claim}

\theoremstyle{definition}

\theoremstyle{remark}
  \newtheorem{rem}[thm]{Remark}
  \newtheorem*{ack}{Acknowledgments}
\newcommand{\Z}{\mathbb{Z}}

\newcommand{\Vol}{\operatorname{Vol}}
\numberwithin{equation}{section}
\date{23rd June, 2002}
\allowdisplaybreaks
\begin{document}
\title[Mahler measure and the volume conjecture]
{Mahler measure of the colored Jones polynomial and the volume conjecture}
\author{Hitoshi Murakami}
\address{
Department of Mathematics,
Tokyo Institute of Technology,
Oh-okayama, Meguro, Tokyo 152-8551, Japan
}
\email{starshea@tky3.3web.ne.jp}
\begin{abstract}
In this note, I will discuss a possible relation between the Mahler measure of
the colored Jones polynomial and the volume conjecture.
In particular, I will study the colored Jones polynomial of the figure-eight
knot on the unit circle.
I will also propose a method to prove the volume conjecture for satellites
of the figure-eight knot.
\end{abstract}
\thanks{This research is partially supported by Grant-in-Aid for Scientific
Research (B)}
\maketitle
\section{Mahler measure}
Let $f(t)$ be a (non-zero) Laurent polynomial in $t$ with coefficient in $\Z$.
The {\em Mahler measure} $\mathbf{M}(f)$ of $f$
\cite{Mahler:MATHE1960,Mahler:JLONM11962,Schmidt:Dynamical_Systems}
is defined to be
\begin{equation*}
  \mathbf{M}(f):=
  \exp
  \left(
    \int_{0}^{1}
    \log
    \left|
      f
      \left(
        \exp(2\pi\sqrt{-1}x)
      \right)
    \right|
    dx
  \right)
\end{equation*}
It is known that $\mathbf{M}(f)$ is the product of the absolute values of the
leading coefficient and all the roots that are greater than one.
It is convenient to define its logarithmic version:
\begin{equation*}
  \mathbf{m}(f):=
    \int_{0}^{1}
    \log
    \left|
      f
      \left(
        \exp(2\pi\sqrt{-1}x)
      \right)
    \right|
    dx.
\end{equation*}
Then the logarithmic Mahler measure can be regarded as a sort of `mean' of the logarithms of the values on the unit circle.
Visit the web pages
\par
http://mathworld.wolfram.com/MahlerMeasure.html
\par\noindent
for more about the Mahler measure
and also
\par
http://math.ucr.edu/~xl/knotprob/knotprob.html
\par\noindent
for problems on the Mahler measure of the Jones polynomial.
\section{Mahler measure of the Alexander polynomial}
Let $K$ be a knot in the three-sphere $S^3$ and $M_N(K)$ be the $N$-fold cyclic
branched covering over $S^3$ branched along $K$.
Then it is well known that the order of the first homology group of $M_N(K)$ can
be obtained in terms of the Alexander polynomial $\Delta(K;t)$ of $K$
(see for example \cite[Corollary~9.8]{Lickorish:1997}).
\begin{thm}
\begin{equation}\label{eq:Nth_roots}
  \left|H_1(M_N(K);\Z)\right|
  =
  \prod_{d=1}^{N-1}\Delta(K,\exp(2d\pi\sqrt{-1}/N)),
\end{equation}
where $|A|$ denotes the cardinality of a set $A$ if $A$ is a finite set and $0$ if it is infinite.
\end{thm}
If we take the logarithm of the both side of Equation~\eqref{eq:Nth_roots}
and divide by $N$, we have
\begin{equation*}
  \frac{\log\left|H_1(M_N(K);\Z)\right|}{N}
  =
  \frac{\sum_{d=1}^{N-1}\log\left|\Delta(K;\exp(2\pi d\sqrt{-1}/N))\right|}{N}
\end{equation*}
When $N$ grows, the right hand side approaches to the `mean' of the logarithms of the values of $\Delta(K;t)$ on the unit circle, the logarithmic Mahler measure.
In fact the following theorem is known to be true.
\begin{thm}[D.~Silver and S.~Williams
\cite{Silver/Williams:Mahler_homology}]
\begin{equation*}
  \lim_{N\to\infty}\frac{\log\left|H_1(M_N(K);\Z)\right|}{N}
  =
  \mathbf{m}(\Delta(K;t))
\end{equation*}
\end{thm}
See \cite{Gordon:TRAAM1972,Gonzalez-Acuna/Short:REVMU1991,Riley:BULLM1990} for other topics of the homology of the branched cyclic cover over a knot.
See also \cite{Silver/Williams:Mahler_Alexander} for the Mahler measure of
the Alexander polynomial of a link.
\section{Mahler measure of the colored Jones polynomials}
Let $J_N(K;t)$ be the $N$-dimensional colored Jones polynomial of a knot $K$
normalized so that $J_N(O;t)=1$ for the unknot $O$.
We want to know the asymptotic behavior of $J_N(K;t)$ for large $N$.
\par
Since
\begin{align*}
  \mathbf{m}\left(J_N(K;t)\right)
  &=
  \int_{0}^{1}\log\left|J_N\bigl(K;\exp(2\pi\sqrt{-1}x)\bigr)\right|dx
  \\
  &=
  \int_{0}^{N}
  \frac{\log\left|J_N\bigl(K;\exp(2r\pi\sqrt{-1}/N)\bigr)\right|}{N}dr,
\end{align*}
it is helpful to study the asymptotic behavior of
$\log\left|J_N\bigl(K;\exp(2\pi\sqrt{-1}r/N)\bigr)\right|$
for a fixed $r$.
Note that for $r=1$, this problem is nothing but the volume conjecture
\cite{Murakami/Murakami:ACTAM12001,
      Murakami:4_1,
      Murakami:Yokota,
      Murakami:optimistic,
      Yokota:Murasugi70,
      Yokota:Topology_Symposium2000,
      Yokota:volume,
      Murakami/Murakami/Okamoto/Takata/Yokota:CS}.
\par
In the following sections I will discuss the colored Jones polynomials of the
figure-eight knot evaluated on the unit circle.
\section{Some calculations about the figure-eight knot}
Let $E$ denote the figure-eight knot $4_1$.
Due to K.~Habiro and T.~Le, the following formula is known.
\begin{equation}\label{eq:fig8}
  J_N(E;t)
  =
  \sum_{k=0}^{N-1}\prod_{j=1}^{k}
  \left(
    t^{(N+j)/2}-t^{-(N+j)/2}
  \right)
  \left(
    t^{(N-j)/2}-t^{-(N-j)/2}
  \right).
\end{equation}
Using this formula we can prove the following result.
\begin{thm}\label{thm}
Let $r$ be a positive integer or a real number satisfying ${5/6}<{r}<{7/6}$.
Then
\begin{equation*}
  \lim_{N\to\infty}
  2\pi\frac{\log\left|J_N\left(E;\exp(2r\pi\sqrt{-1})\right)\right|}{N}
  =
  \frac{2\Lambda\bigl(r\pi+\theta(r)/2\bigr)
       -2\Lambda\bigl(r\pi-\theta(r)/2\bigr)}{r},
\end{equation*}
where
$\Lambda(z):=-\displaystyle\int_{0}^{z}\log|\sin{x}|dx$
is the Lobachevski function
and $\theta(r)$ is the smallest positive number satisfying
$\cos\theta(r)=\cos(2r\pi)-1/2$.
\par
In particular, if $r$ is a positive integer then
\begin{equation*}
  2\pi\lim_{N\to\infty}
  \frac{\log\left|J_N\left(E;\exp(2r\pi\sqrt{-1}/N)\right)\right|}{N}
  =
  \frac{\Vol(S^3\setminus{E})}{r}.
\end{equation*}
\end{thm}
\begin{proof}[Proof of Theorem~\ref{thm} when $r$ is a positive integer]
\par
Replacing $t$ with $\exp(2r\pi\sqrt{-1}/N)$ in Equation~\eqref{eq:fig8}, we have
\begin{equation*}
  J_N(E;\exp(2r\pi\sqrt{-1}/N))
  =
  \sum_{k=0}^{N-1}\prod_{j=1}^{k}
  \left\{2\sin(jr\pi/N)\right\}^2
\end{equation*}
If we put $f(k):=\prod_{j=1}^{k}\left\{2\sin(jr\pi/N)\right\}^2$,
then $f$ takes its maximum at $kr\pi/N=5\pi/6$ if $N$ is large.
Therefore
\begin{align*}
  \lim_{N\to\infty}
  \frac{\log\left|J_N\left(E;\exp(2r\pi\sqrt{-1}/N)\right)\right|}{N}
  &=
  2\lim_{N\to\infty}\sum_{j=1}^{5N/6r}
  \frac{\log\left(2\sin(jr\pi/N)\right)}{N}
  \\
  &=
  \frac{2}{r\pi}\lim_{N\to\infty}
  \int_{0}^{5\pi/6}\log\left(2\sin{x}\right)\;dx
  \\
  &=
  -\frac{2}{r\pi}\Lambda(5\pi/6)
  \\
  &=
  \frac{\Vol(S^3\setminus{E})}{2r\pi}.
\end{align*}
See \cite[Theorem 4.2]{Murakami:4_1} for details.
\end{proof}
\begin{rem}
The case where $r=1$ is due to R.~Kashaev \cite{Kashaev:LETMP97}
and T.~Ekholm \cite{Murakami:4_1}.
\end{rem}
\begin{proof}[Proof of Theorem~\ref{thm} when $5/6<r<1$]
We will assume $N$ is sufficiently large so that $j/N$ can behave as if it is a
continuous parameter.
\par
Put $\omega:=\exp(2\pi\sqrt{-1}/N)$.
Since
\begin{align*}
  \omega^{r(N+j)/2}-\omega^{-r(N+j)/2}&=2\sqrt{-1}\sin(r(N+j)\pi/N)
  \\
  \intertext{and}
  \omega^{r(N-j)/2}-\omega^{-r(N-j)/2}&=2\sqrt{-1}\sin(r(N-j)\pi/N),
\end{align*}
we have
\begin{multline*}
  \prod_{j=1}^{k}
  \left(
    \omega^{r(N+j)/2}-\omega^{-r(N+j)/2}
  \right)
  \left(
    \omega^{r(N-j)/2}-\omega^{-r(N-j)/2}
  \right)
  \\
  =
  \prod_{j=1}^{k}
  4\sin(rj\pi/N+r\pi)\sin(rj\pi/N-r\pi).
\end{multline*}
Put
\begin{align*}
g(j):&=4\sin(rj\pi/N+r\pi)\sin(rj\pi/N-r\pi)
\\
     &=2\cos(2r\pi)-2\cos(2rj\pi/N)
\\
\intertext{and}
f(k):&=\prod_{j=1}^{k}g(j)
\end{align*}
so that
$J_N(E;\omega^r)=\sum_{k=0}^{N-1}f(k)$.
\begin{figure}[h]
\begin{center}
\includegraphics[scale=0.45]{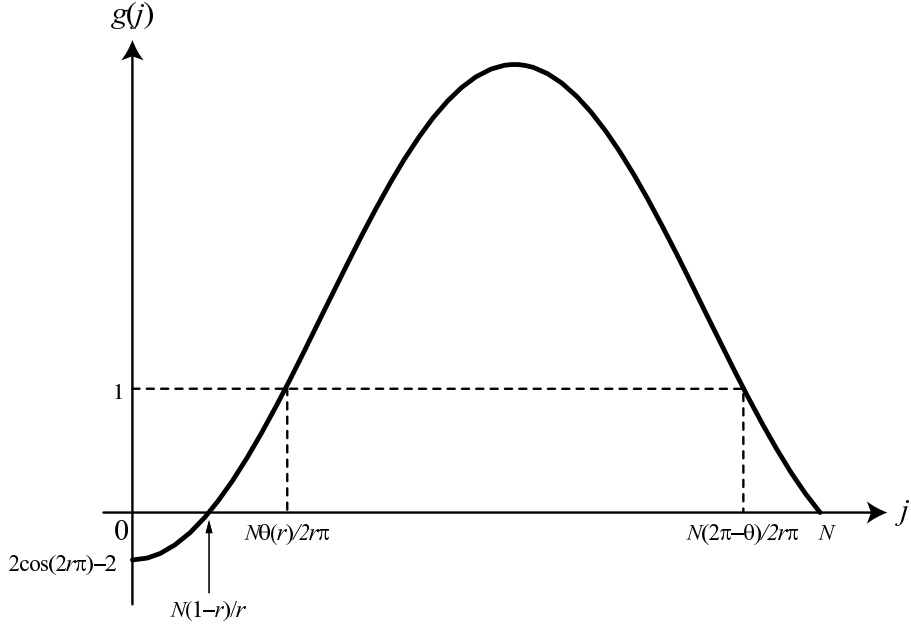}
\caption{Graph of $g(j)$ when $5/6<r<1$}
\label{fig:cos1}
\end{center}
\end{figure}
We also put
\begin{equation*}
  A:=\frac{N(1-r)}{r},\quad
  B:=\frac{N\theta(r)}{2r\pi},\quad
  \text{and}\qquad
  C:=\frac{N(2\pi-\theta(r))}{2r\pi}
\end{equation*}
where $\theta(r)$ is the smallest positive
number satisfying $\cos\theta(r)=\cos(2r\pi)-1/2$ as before.
Note that since $5/6<r<1$, $1/2<\cos(2r\pi)<1$ and so the equation
$\cos\theta(r)=\cos(2r\pi)-1/2$
has a solution.
\par
Note that $0<A<B<C<N$ (see Figure~\ref{fig:cos1}). 
\par
Since we have
\begin{enumerate}
\item
  $g(j)<0$ for $j<A$, and $g(j)>0$ for $j>A$,
\quad and
\item
  $f_j>1$ for $B<j<C$,
\end{enumerate}
we see
\begin{enumerate}
\item[(3)]
  If $j<A$ then the signs of $f(j)$ alternate, that is,
  $f(j-1)f(j)<0$, and if $j>A$ then the signs of $f(j)$ are constant,
\quad and
\item[(4)]
  $|f(0)|>|f(1)|>\dots>|f(B)|$ and $|f(B+1)|<\dots<|f(C)|$.
\end{enumerate}
\par
Let $f_{\text{MAX}}$
\footnote{MAX are temporarily Nana, Reina, and Lina.}
be the maximum of $\{|f_j|\}$ for $0\le j \le N-1$.
Note that $f_{\text{MAX}}=f(C)$.
We can show the following inequality.
\begin{claim}\label{claim:inequality1}
\begin{equation*}
  0<f_{\text{MAX}}-1 \le \left|J_N(E;\omega^r)\right| \le N\,f_{\text{MAX}}
\end{equation*}
\end{claim}
\begin{proof}[Proof of the Claim~\ref{claim:inequality1}]
We only show the second inequality for the case where $A$ is even.
In this case since $f(0)=1$, $f(2j-1)+f(2j)<0$ for $2j<A$, and $f(j)<0$ for
$j\ge{A-1}$, we have
\begin{align*}
  &\left|J_N(E;\omega^r)\right|
  \\
  &\quad=
  \bigl|f(0)+
  \{f(1)+f(2)\}+\{f(3)+f(4)\}+\dots+\{f(A-3)+f(A-2)\}
  \\
  &\quad\quad
  +f(A-1)+f(A)+f(A+1)+\dots+f(N-1)
  \bigr|
  \\
  &\quad=
  |f(1)+f(2)|+|f(3)+f(4)|+\dots+|f(A-3)+f(A-2)|
  \\
  &\quad\quad
  +|f(A-1)|+|f(A)|+\dots+|f(N-1)|
  \\
  &\quad\quad
  -1
  \\
  &\quad>
  f_{\text{MAX}}-1
\end{align*}
and the second equality follows.
\end{proof}
Therefore we have
\begin{align*}
  &\lim_{N\to\infty}\frac{\log|J_N(E;\omega^r)|}{N}
  \\
  &\quad=
  \lim_{N\to\infty}\frac{\log(f_{\text{MAX}})}{N}
  \\
  &\quad=
  \lim_{N\to\infty}\frac{1}{N}
  \sum_{j=0}^{C}
  \left\{
     \log\left(2\sin\left(\frac{rj\pi}{N}+r\pi-\pi\right)\right)
    +\log\left(2\sin\left(-\frac{rj\pi}{N}+r\pi\right)\right)
  \right\}
  \\
  &\quad=
   \frac{1}{r\pi}
   \int_{r\pi-\pi}^{r\pi-\theta(r)/2}\log(2\sin{x})d\,x
  +\frac{1}{r\pi}
   \int_{r\pi-\pi+\theta(r)/2}^{r\pi}\log(2\sin{x})d\,x
  \\
  &\quad=
  \frac{1}{r\pi}
  \left(
     \Lambda(r\pi-\pi)-\Lambda(r\pi-\theta(r)/2)
    +\Lambda(r\pi-\pi+\theta(r)/2)-\Lambda(r\pi)
  \right)
  \\
  &\quad=
  \frac{1}{r\pi}
  \left(
    \Lambda(r\pi+\theta(r)/2)-\Lambda(r\pi-\theta(r)/2)
  \right).
\end{align*}
Here we use the $\pi$-periodicity of the Lobachevski function.
(See \cite{Milnor:BULAM382}.)
\end{proof}
\begin{proof}[Proof when $1<r<7/6$]
The proof is similar to the case where $5/6<r<1$.
See Figure~\ref{fig:cos2}.
\begin{figure}[h]
\begin{center}
\includegraphics[scale=0.45]{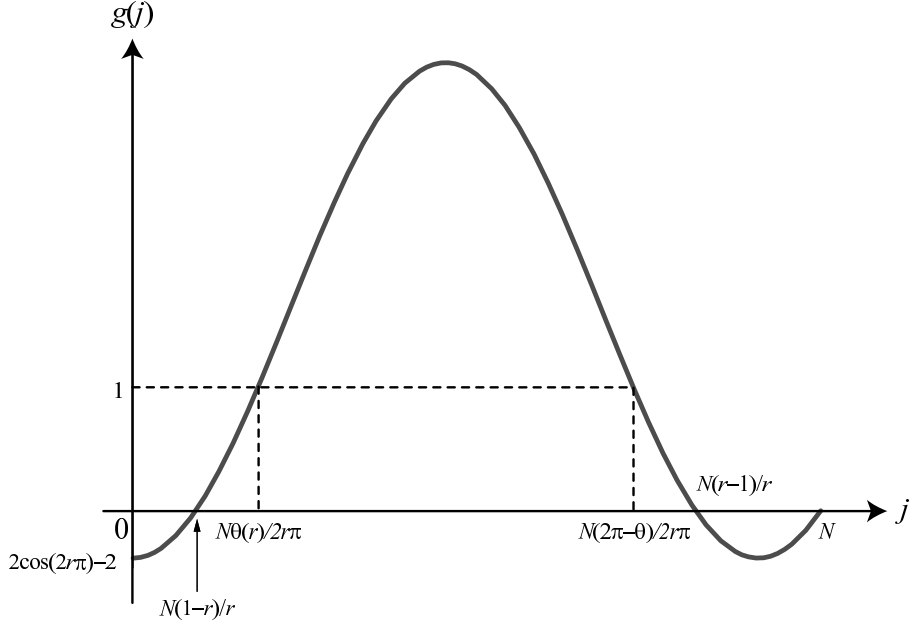}
\caption{Graph of $g(j)$ when $1<r<7/6$}
\label{fig:cos2}
\end{center}
\end{figure}
\end{proof}
As a corollary we have
\begin{cor}\label{cor:limit}
\begin{equation*}
  \lim_{r\to1}
  \left\{
    \lim_{N\to\infty}
    \frac{\log\left|J_N(E;\exp(2r\pi\sqrt{-1}))\right|}{N}
  \right\}
  =
    \lim_{N\to\infty}
    \frac{\log\left|J_N(E;\exp(2\pi\sqrt{-1}))\right|}{N}
\end{equation*}
\end{cor}
By some calculation using PARI-GP
\footnote{\par\vspace{-3mm}
\begin{minipage}{10cm}
\begin{center}
\texttt{GP/PARI CALCULATOR Version 2.0.20 (beta)}
\par
\texttt{i586 running cygwin\_98-4.10 (ix86 kernel) 32-bit version}
\par
\texttt{(readline v1.0 enabled, extended help not available)}
\par\quad\par
\texttt{Copyright (C) by 1989-1999 by}
\par
\texttt{C. Batut, K. Belabas, D. Bernardi, H. Cohen and M. Olivier.}
\end{center}
\par\medskip
The program is available at
http://www.parigp-home.de/
\end{minipage}
}
and MAPLE V, it seems that the following equality holds.
\begin{equation}\label{eq:V_W}
  2r\pi\lim_{N\to\infty}\frac{\log\left|J_N(E;\omega^r)\right|}{N}
  =
  \begin{cases}
    V(r)    \quad&\text{if ${0}\le{r}\le{1}$},
    \\
    W(r-[r])\quad&\text{if $r>1$},
  \end{cases}
\end{equation}
where $[r]$ denotes the greatest integer which does not exceed $r$,
and
\begin{align*}
  V(x):=
  \begin{cases}
    0
    \quad&\text{if ${0}\le{x}<{1/6}$,}
    \\[2mm]
    \Lambda(x\pi+\theta(x)/2-\pi/2)-\Lambda(x\pi-\theta(x)/2-\pi/2)
    \quad&\text{if ${1/6}\le{x}<3/4$,}
    \\[2mm]
    \Lambda(x\pi+\theta(x)/2)-\Lambda(x\pi-\theta(x)/2)
    \quad&\text{if ${3/4}\le{x}\le{1}$,}
  \end{cases}
  \\
  \intertext{and}
  W(x):=
  \begin{cases}
    \Lambda(x\pi)+\theta(x)/2-\Lambda(x\pi-\theta(x)/2)
    \quad&\text{if ${0}\le{x}<{1/4}$,}
    \\[2mm]
    \Lambda(x\pi)+\theta(x)/2-\pi/2)-\Lambda(x\pi-\theta(x)/2-\pi/2)
    \quad&\text{if ${1/4}\le{x}<3/4$,}
    \\[2mm]
    \Lambda(x\pi+\theta(x)/2)-\Lambda(x\pi-\theta(x)/2)
    \quad&\text{if ${3/4}\le{x}\le{1}$.}
  \end{cases}
\end{align*}
See Figures~\ref{fig:V} and \ref{fig:W} for graphs of $V$ and $W$.
\begin{figure}[h]
\begin{center}
\includegraphics[scale=0.45]{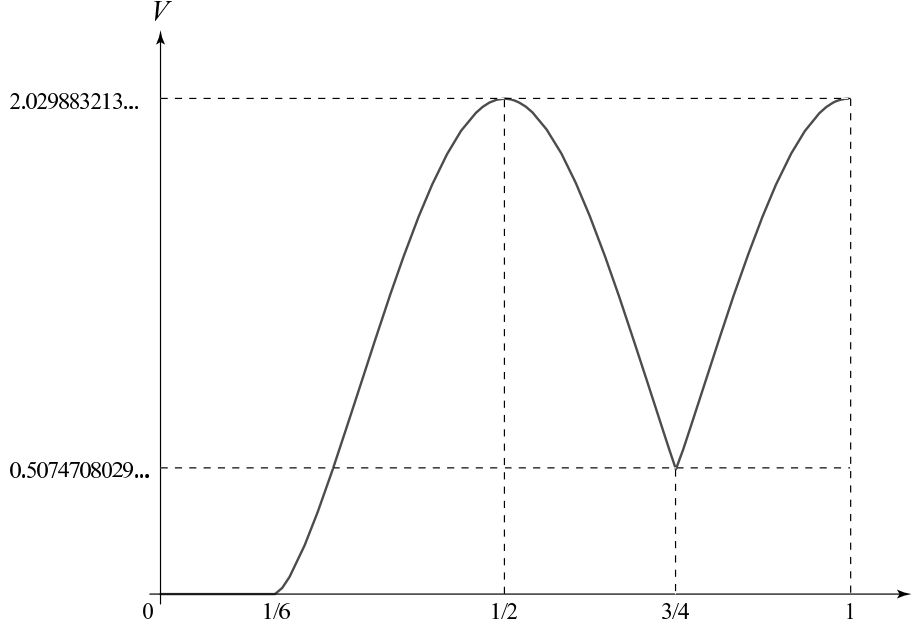}
\caption{Graph of $V$, where $2.029883213...$ is the volume of the figure-eight
knot complement.}
\label{fig:V}
\end{center}
\end{figure}
\begin{figure}[h]
\begin{center}
\includegraphics[scale=0.45]{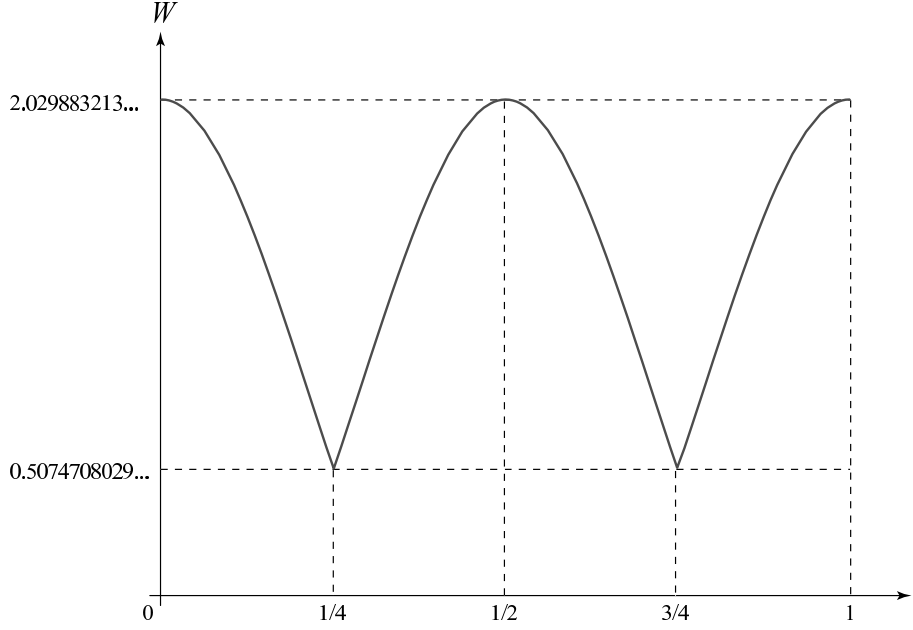}
\caption{Graph of $W$.}
\label{fig:W}
\end{center}
\end{figure}
See also
Figures~\ref{fig:2000_1000_1},
        \ref{fig:2000_1000_2},
        \ref{fig:2000_1000_3},
        \ref{fig:2000_1000_4},
        \ref{fig:2000_1000_5}, and
        \ref{fig:8000_1000_5}
for some results of calculations supporting
Equation~\ref{eq:V_W}.
\par
If Equation~\eqref{eq:V_W} is true, one could have the following result on
the asymptotic behavior of the logarithmic Mahler measure of the colored
Jones polynomials of the figure-eight knot.
\begin{rem}
Caution!
There are {\em fake} calculations in the following.
\end{rem}
\begin{align*}
  \lim_{N\to\infty}\frac{\mathbf{m}\left(J_N(E;t)\right)}{\log{N}}
  &=
  \lim_{N\to\infty}\frac{1}{\log{N}}
  \int_{0}^{1}\log\left|J_N(E;\exp(2\pi\sqrt{-1}x))\right|dx
  \\
  &\underset{?}{=}
  \lim_{N\to\infty}\frac{1}{\log{N}}
  \int_{0}^{N}
  \frac{\log\left|J_N(E;\exp(2\pi\sqrt{-1}r/N))\right|}{N}dr
  \\
  &\underset{?}{=}
  \lim_{N\to\infty}\frac{1}{2\pi\log{N}}
  \left\{
    \int_{0}^{1}\frac{V(r)}{r}dr
    +
    \sum_{k=1}^{N-1}\int_{k}^{k+1}\frac{W(r-[r])}{r}dr
  \right\}
  \\
  &=
  \lim_{N\to\infty}\frac{1}{2\pi\log{N}}
  \left\{
    \int_{0}^{1}\frac{V(r)}{r}dr
    +
    \sum_{k=1}^{N-1}\int_{0}^{1}\frac{W(r)}{r+k}dr
  \right\},
\end{align*}
where $\underset{?}{=}$ means that there is a doubt in the equality.
At the first I use $N$ in the integral, which should be independent of $N$,
and at the second I assume \eqref{eq:V_W}.
\par
Now since
\begin{equation*}
  \frac{1}{k+1}\le\frac{1}{r+k}\le\frac{1}{k}
\end{equation*}
for ${0}\le{r}\le{1}$, we have
\begin{equation*}
  \int_{0}^{1}\frac{W(r)}{k+1}dr
  \le
  \int_{0}^{1}\frac{W(r)}{r+k}dr
  \le
  \int_{0}^{1}\frac{W(r)}{k}dr.
\end{equation*}
Therefore we have
\begin{equation*}
  \sum_{k=1}^{N-1}
  \int_{0}^{1}\frac{W(r)}{k+1}dr
  \le
  \sum_{k=1}^{N-1}
  \int_{0}^{1}\frac{W(r)}{r+k}dr
  \le
  \sum_{k=1}^{N-1}
  \int_{0}^{1}\frac{W(r)}{k}dr.
\end{equation*}
Since
\begin{equation*}
  \lim_{N\to\infty}\frac{\sum_{k=1}^{N-1}\frac{1}{k+1}}{\log{N}}
  =
  \lim_{N\to\infty}\frac{\sum_{k=1}^{N-1}\frac{1}{k  }}{\log{N}}
  =
  1
\end{equation*}
and $V(r)=0$ for ${0}\le{r}\le{1/6}$, we finally have
\begin{equation*}
  2\pi\lim_{N\to\infty}\frac{\mathbf{m}\left(J_N(E;t)\right)}{\log{N}}
  \underset{?}{=}
  \int_{0}^{1}W(r)dr
  =1.450191516....
\end{equation*}
\section{Satellites of the figure-eight knot}
In this section, I would like to study the volume conjecture for the
$(2,1)$-cable and the Whitehead double of the figure-eight knot.
Linear skein method gives us formulas to describe the colored Jones polynomials
of such knots but one of the difficulties is that the value of the unknot
is not $1$ but $(t^{N/2}-t^{-N/2})/(t^{1/2}-t^{-1/2})$
(see for example \cite[Chapter~14]{Lickorish:1997}), and so they vanish
if we evaluate them at the $N$-th root of unity.
To avoid this I will use Corollary~\ref{cor:limit} to analyze the asymptotic
behaviors of the colored Jones polynomials.
Unfortunately, I cannot give a rigorous result here but I hope that this method
gives an insight to solve the volume conjecture for satellite knots.
\begin{rem}
Caution!
There are many {\em fake} arguments in this section.
\end{rem}
\par
\begin{figure}[h]
\begin{center}
\includegraphics[scale=0.6]{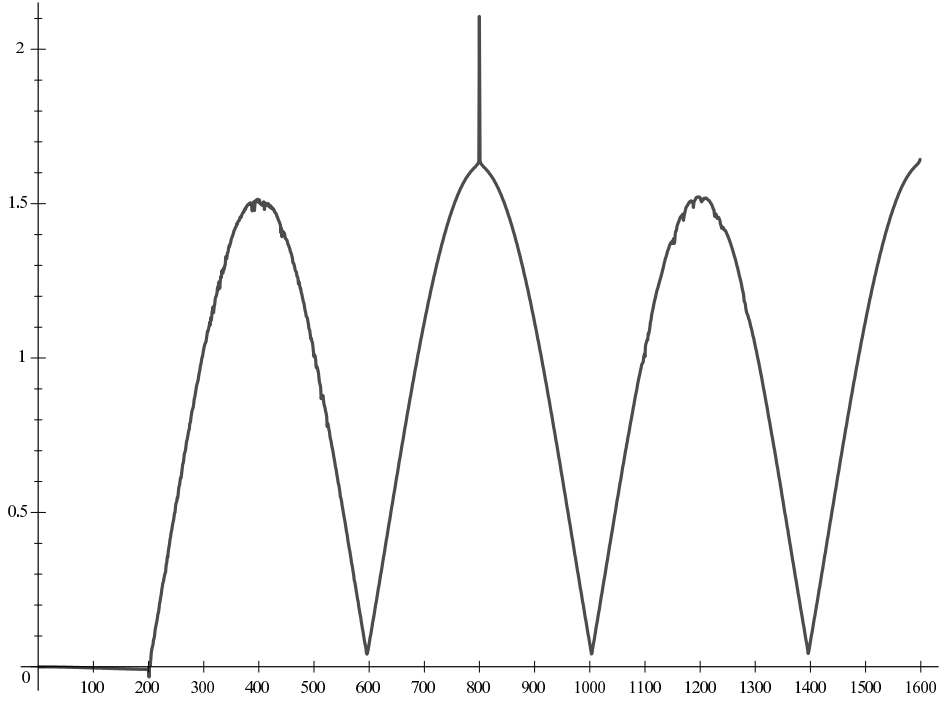}
\caption{Plot of $(c,2\pi\log\left|J_c(E;\omega)\right|/N)$ with $N=800$.}
\label{fig:max}
\end{center}
\end{figure}
\begin{figure}[h]
\begin{center}
\includegraphics[scale=0.45]{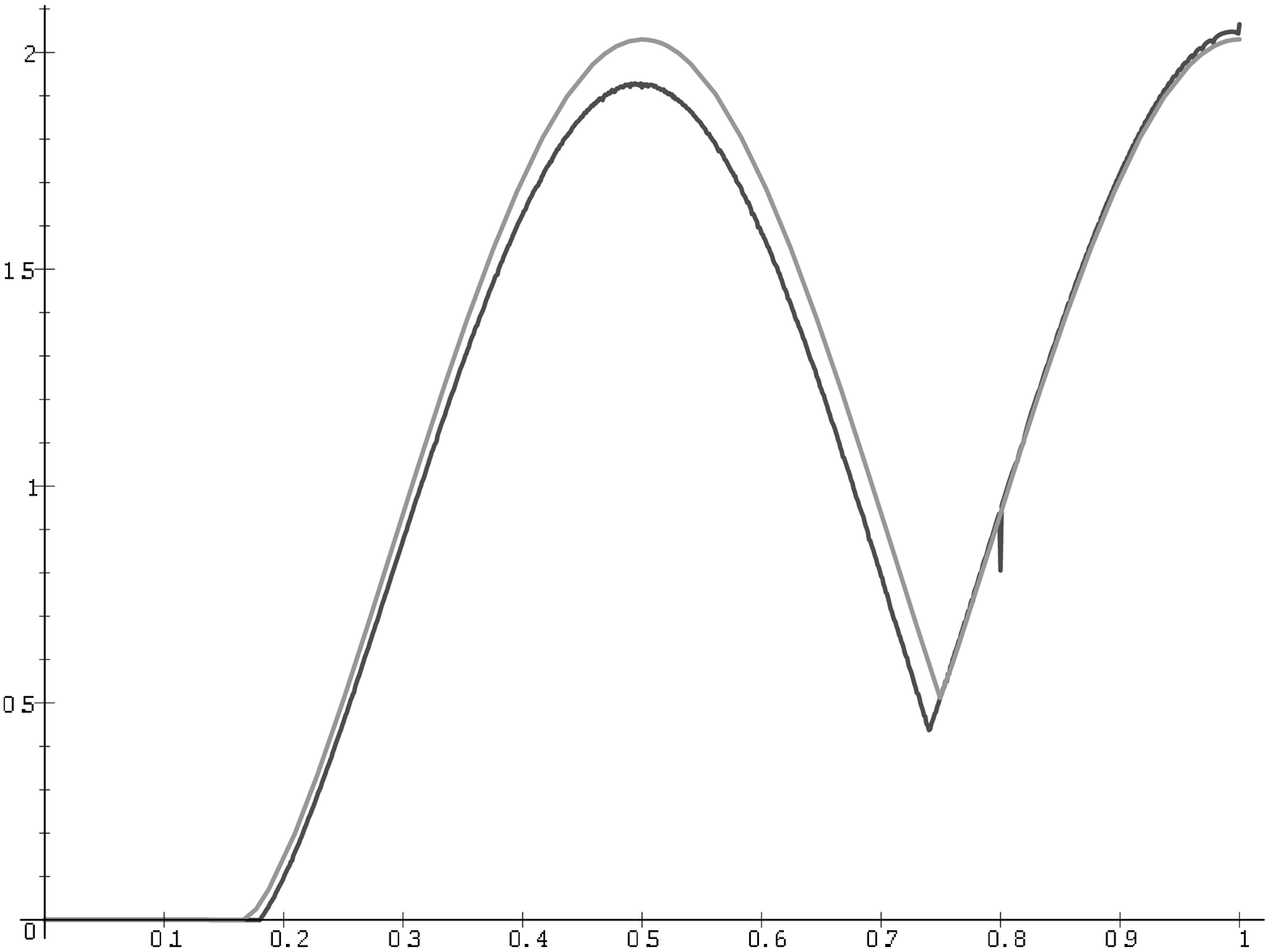}
\caption{Graph of $W$ (gray) and $2r\pi\log\left|J_N(E;\omega^r)\right|/N$
with $N=2000$ (black) for $0\le{r}\le1$.}
\label{fig:2000_1000_1}
\end{center}
\end{figure}
\begin{figure}[h]
\begin{center}
\includegraphics[scale=0.45]{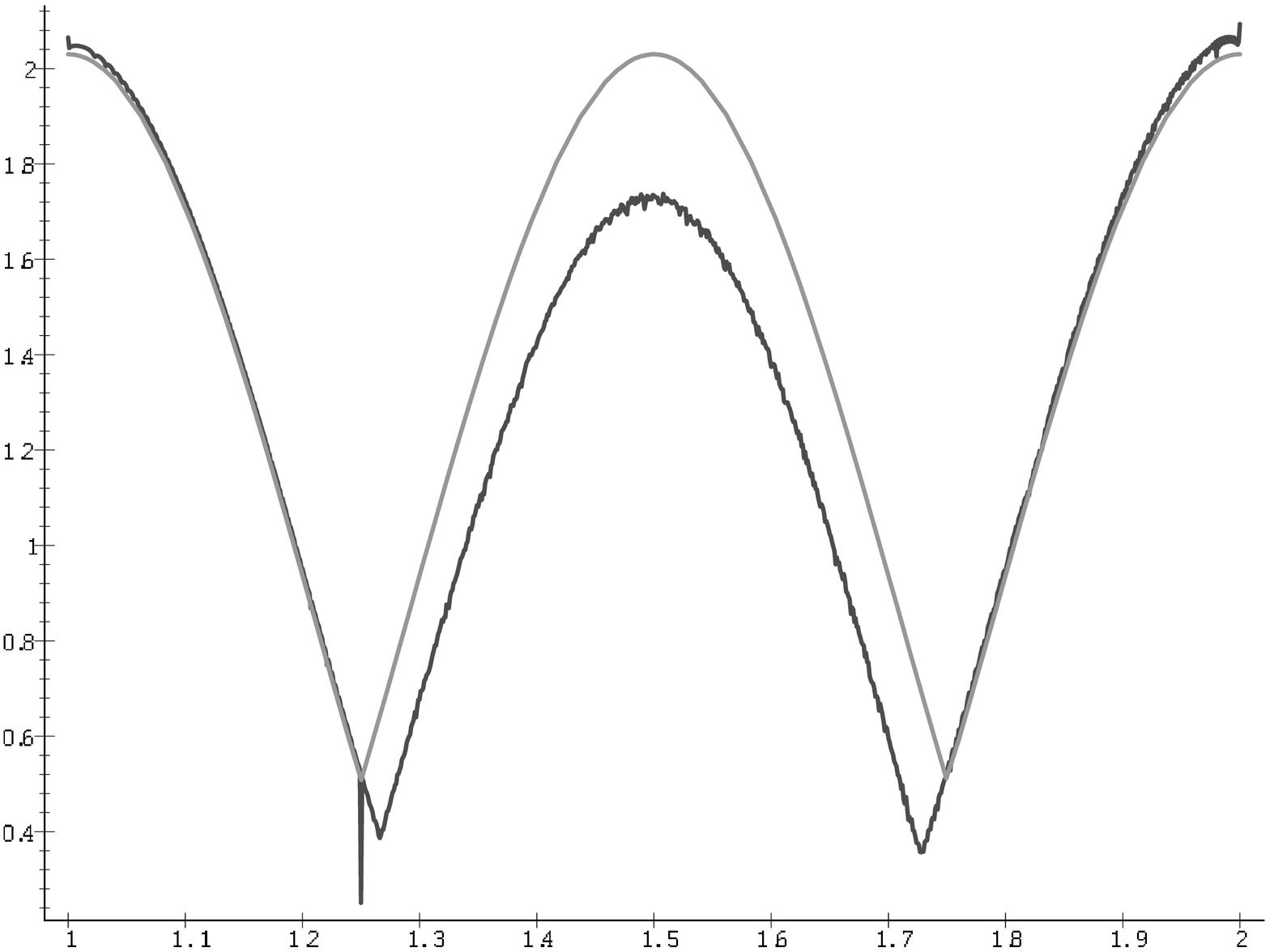}
\caption{Graph of $W$ (gray) and $2r\pi\log\left|J_N(E;\omega^r)\right|/N$
with $N=2000$ (black) for $1\le{r}\le2$.}
\label{fig:2000_1000_2}
\end{center}
\end{figure}
\begin{figure}[h]
\begin{center}
\includegraphics[scale=0.45]{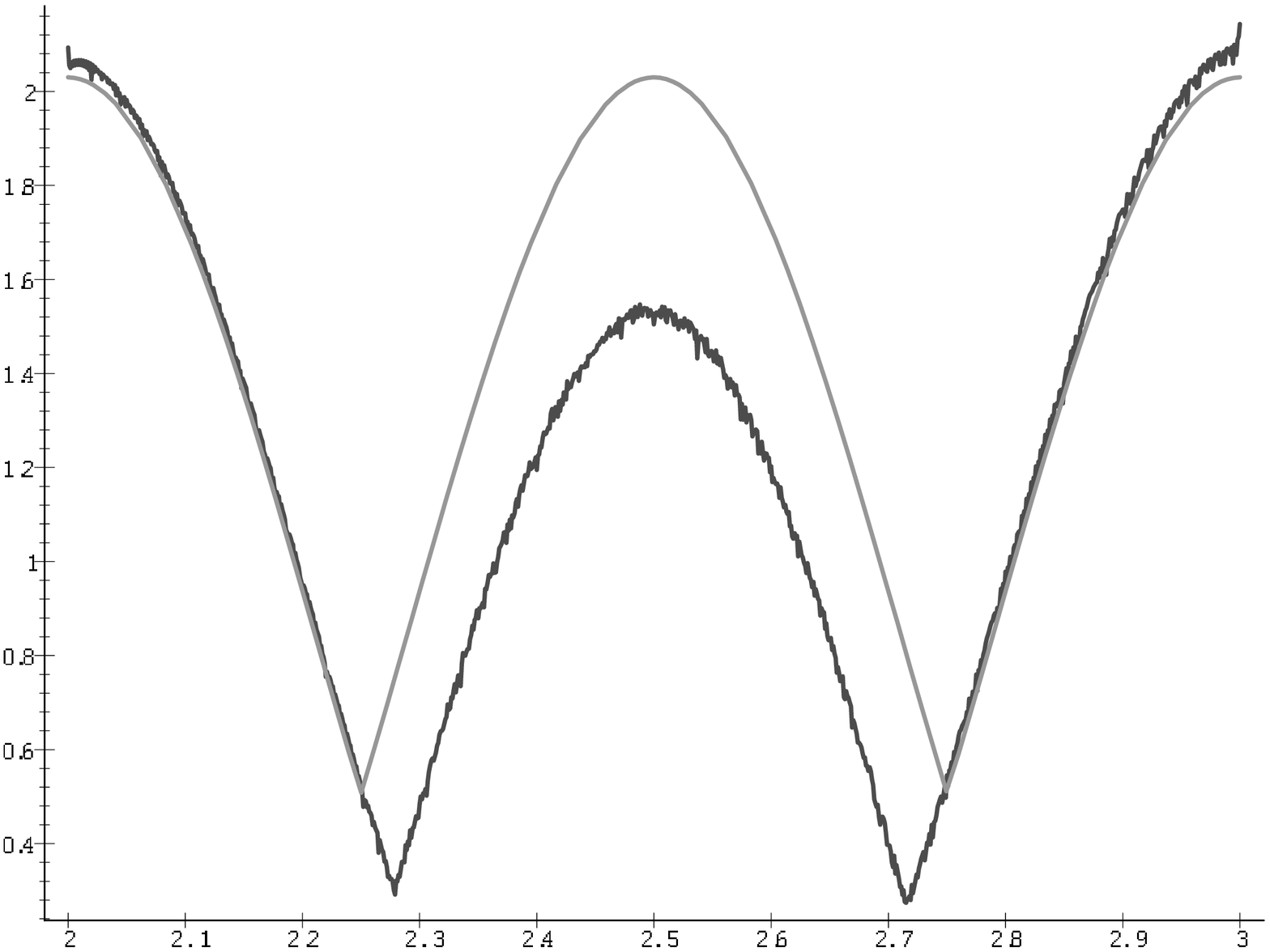}
\caption{Graph of $W$ (gray) and $2r\pi\log\left|J_N(E;\omega^r)\right|/N$
with $N=2000$ (black) for $2\le{r}\le3$.}
\label{fig:2000_1000_3}
\end{center}
\end{figure}
\begin{figure}[h]
\begin{center}
\includegraphics[scale=0.45]{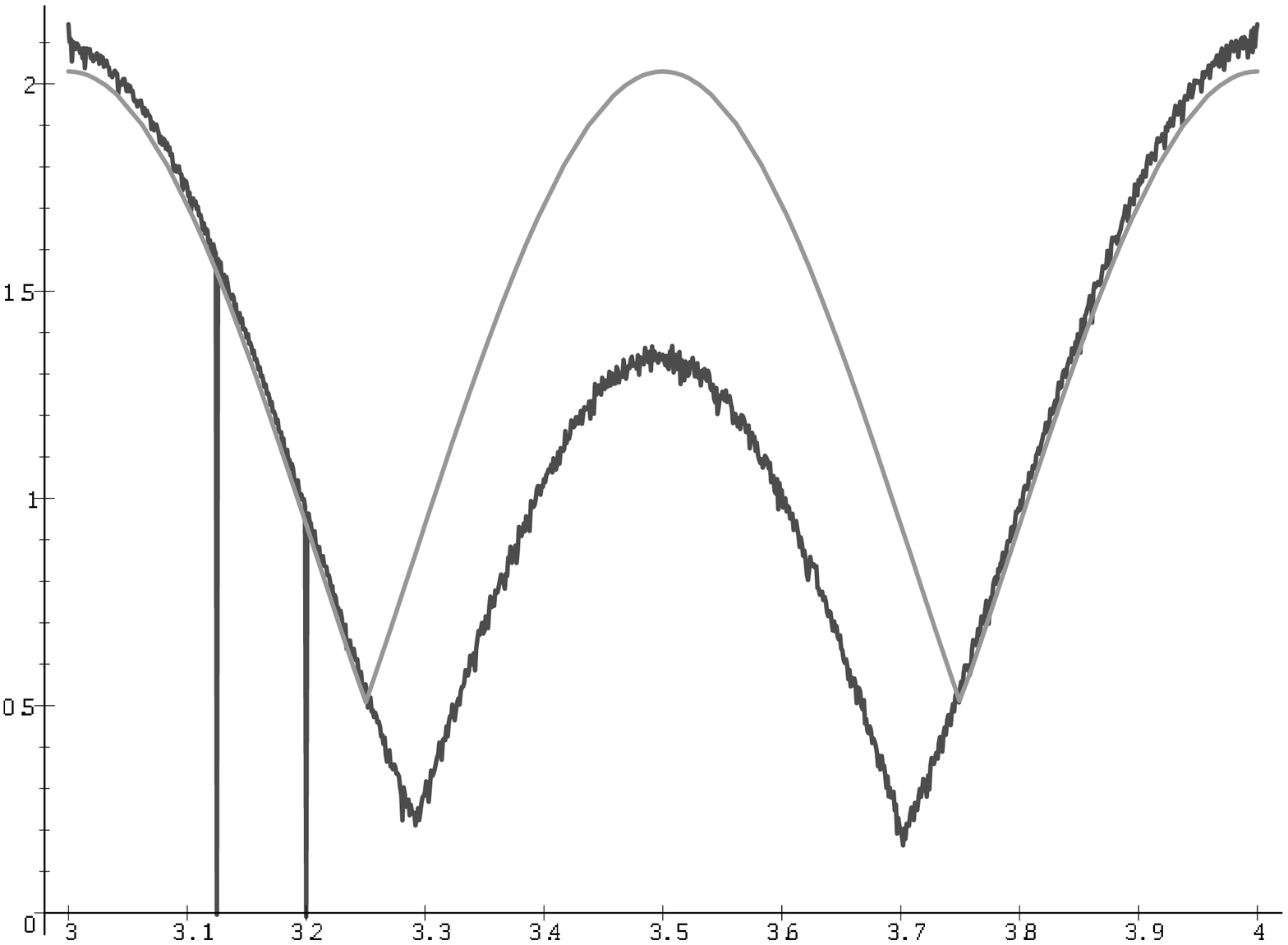}
\caption{Graph of $W$ (gray) and $2r\pi\log\left|J_N(E;\omega^r)\right|/N$
with $N=2000$ (black) for $3\le{r}\le4$.}
\label{fig:2000_1000_4}
\end{center}
\end{figure}
\begin{figure}[h]
\begin{center}
\includegraphics[scale=0.45]{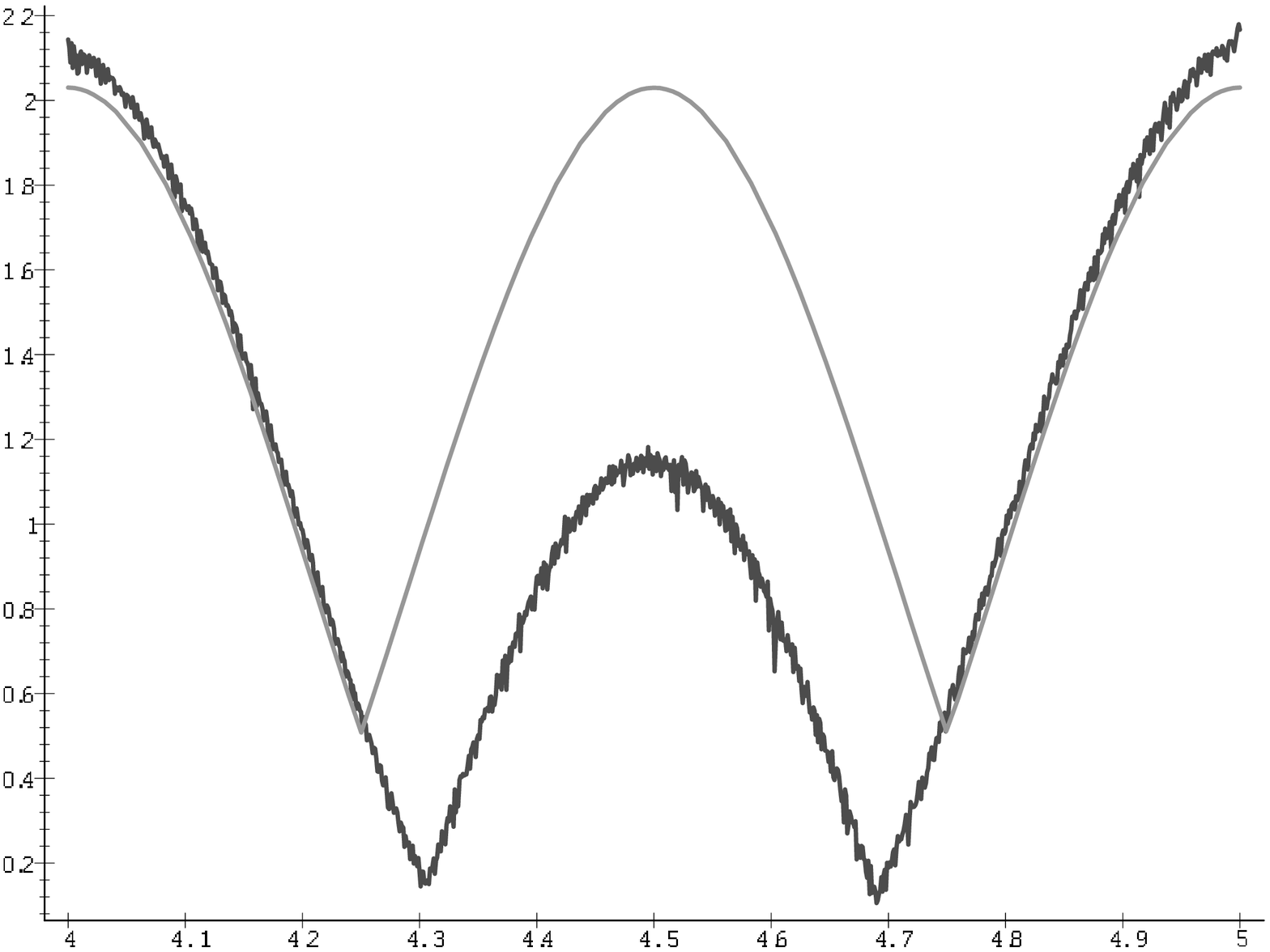}
\caption{Graph of $W$ (gray) and $2r\pi\log\left|J_N(E;\omega^r)\right|/N$
with $N=2000$ (black) for $4\le{r}\le5$.}
\label{fig:2000_1000_5}
\end{center}
\end{figure}
\begin{figure}[h]
\begin{center}
\includegraphics[scale=0.45]{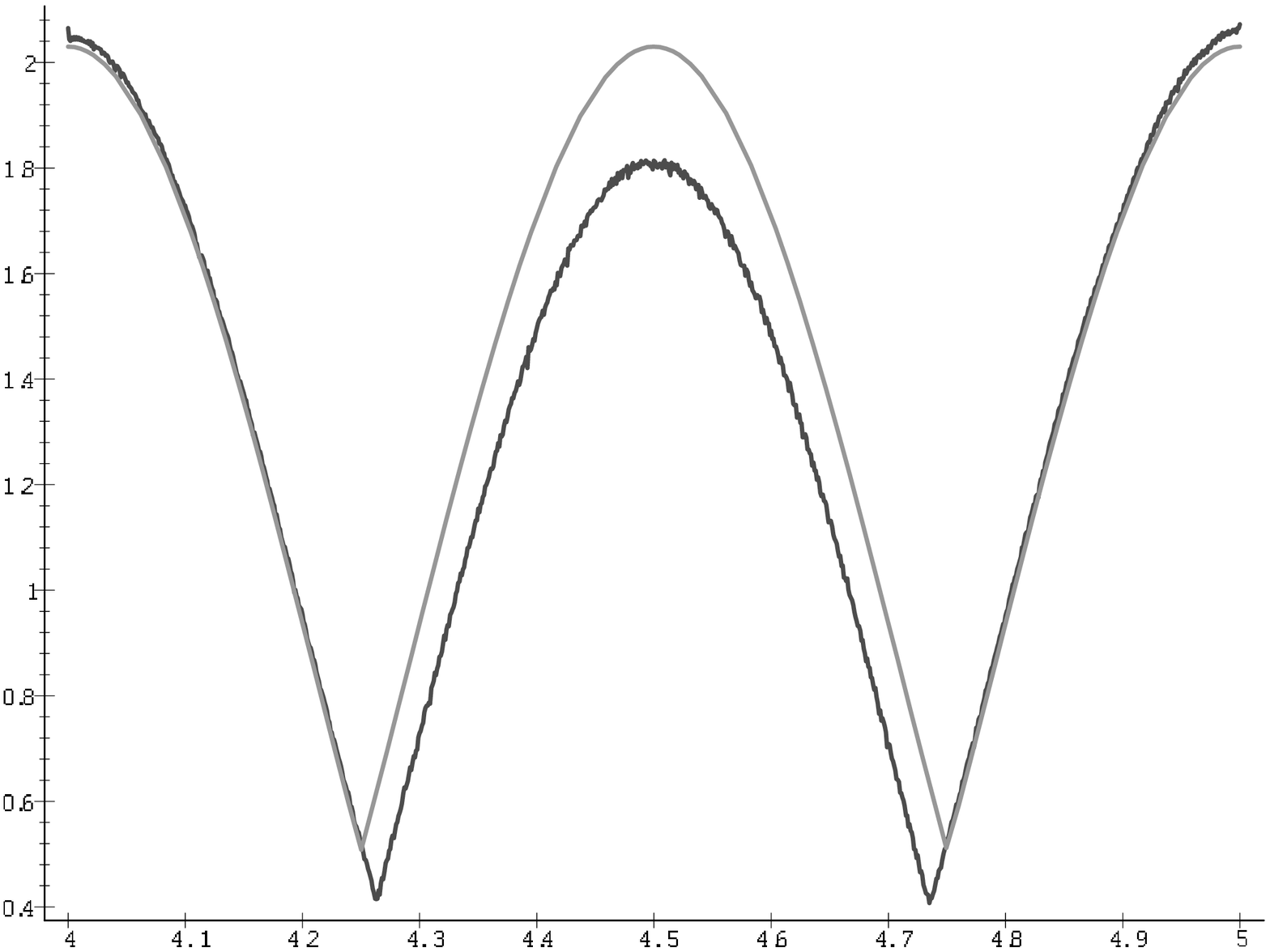}
\caption{Graph of $W$ (gray) and $2r\pi\log\left|J_N(E;\omega^r)\right|/N$
with $N=8000$ (black) for $4\le{r}\le5$.}
\label{fig:8000_1000_5}
\end{center}
\end{figure}

Let $E^2$ be the $(2,1)$-cable of the figure-eight knot.
By using techniques in \cite[Chapter~14]{Lickorish:1997}, we see
\begin{equation*}
  J_N(E^2;t)(t^{N/2}-t^{-N/2})/(t^{1/2}-t^{-1/2})
  =
  \sum_{\substack{c:\text{ odd} \\ c\le 2N-1}}
  u(c;t^{1/4})
  J_c(E;t),
\end{equation*}
where $u(c;t^{1/4})$ is a monomial in $t^{1/4}$.
Replacing $t$ with $\omega^r$ with $5/6<r<7/6$ ($r\ne1$), we have
\begin{equation*}
  J_N(E^2;\omega^r)
  =
  \frac{\sin(r\pi/N)}{\sin(r\pi)}
  \sum_{\substack{c:\text{ odd} \\ 1\le{c}\le 2N-1}}
  u(c;\omega^{r/4})
  J_c(E;\omega^{r}).
\end{equation*}
Note that $\sin(r\pi)\ne0$.
If one could show that the maximum of the terms in the summation dominates the
limit, which is a kind of saddle point method, we
could have
\begin{align*}
  \lim_{N\to\infty}
  \frac{\log\left|J_N(E^2;\omega)\right|}{N}
  &\underset{?}{=}
  \lim_{r\to1}
  \left\{
    \lim_{N\to\infty}
    \frac{\log\left|J_N(E^2;\omega^r)\right|}{N}
  \right\}
  \\
  &\underset{?}{=}
  \lim_{r\to1}
  \left\{
    \frac{\log\left|
                \displaystyle\max_{1\le c\le 2N-1}{J_c(E,\omega^{r})}
              \right|}{N}
  \right\}
  \\
  &\underset{?}{=}
  \lim_{r\to1}
  \left\{
    \frac{\log\left|J_N(E,\omega^{r})\right|}{N}
  \right\}
  \\
  &=
    \frac{\log\left|J_N(E,\omega)\right|}{N},
\end{align*}
proving the volume conjecture for the $(2,1)$-cable of the figure-eight knot.
Here $\underset{?}{=}$ indicates that there is a doubt in the equality;
at the fist equality, I change the order of the limits, 
at the second, I assume the maximum dominates the limit,
and at the third, I assume
that $J_c(E,\omega^r)$ takes its maximum at $c=N$, which can be observed
by calculation using PARI.
See Figure~\ref{fig:max}.
\par
I believe that the gaps here are not so big.
\par
Let $D(E)$ be the Whitehead double of the figure-eight knot (with any framing).
Then using similar techniques we have
\begin{equation*}
  J_N(D(E);t)(t^{N/2}-t^{-N/2})/(t^{1/2}-t^{-1/2})
  =
  \sum_{\substack{c:\text{ odd} \\ c\le 2N-1}}
  v(c;t)J_{c}(E;t),
\end{equation*}
with
\begin{equation*}
  v(c;t)
  =
  \sum_{\substack{d:\text{ odd} \\ d \le 2N-1}}
  \frac{\Delta_{d}\;\theta(N-1,N-1,c-1)}
       {\Delta_{c}\;\theta(N-1,N-1,d-1)}
  \left\{
    \begin{matrix} N-1 & N-1 & c \\ N-1 & N-1 & d\end{matrix}
  \right\}
\end{equation*}
where $\Delta_{x},\theta(x,y,z)$ and
$\left\{\begin{matrix}x&y&z\\u&v&w\end{matrix}\right\}$
are defined in \cite[Chapter~14]{Lickorish:1997}.
Similar calculation shows that for the Whitehead link $W$ we have
\begin{equation*}
  J_N(W;t)(t^{N/2}-t^{-N/2})/(t^{1/2}-t^{-1/2})
  =
  \sum_{\substack{c:\text{ odd} \\ c \le 2N-1}}
  v(c;t)J_{N,c}(H;t),
\end{equation*}
where $J_{N,c}(H;t)$ is the colored Jones polynomial of the Hopf link $H$
colored with $N$ and $c$, which is equal to $\Delta_{(N-1)(c-1)}$.
\par
Now we have the following {\em fake} calculations with doubtful equalities:
\begin{align*}
  \lim_{N\to\infty}
  \frac{\log\left|J_N(D(E);\omega)\right|}{N}
  &\underset{?}{=}
  \lim_{r\to1}
  \left\{
    \lim_{N\to\infty}
    \frac{\log\left|J_N(D(E);\omega^r)\right|}{N}
  \right\}
  \\
  &\underset{?}{=}
  \lim_{r\to1}
  \left\{
    \frac{\log\left|
                \displaystyle\max_{1\le c\le 2N-1}
                             {v(c;\omega^r)J_c(E,\omega^{r})}
              \right|}{N}
  \right\}
  \\
  &\underset{?}{=}
  \lim_{r\to1}
  \left\{
    \frac{\log\left|v(N;\omega^r)J_{N}(E;\omega^r)\right|}{N}
  \right\}
  \\
  &=
  \lim_{r\to1}\frac{\log\left|v(N;\omega^r)\right|}{N}
  +
  \lim_{r\to1}\frac{\log\left|J_N(e;\omega^r)\right|}{N}
  \\
  &=
  \lim_{r\to1}\frac{\log\left|v(N;\omega^r)\right|}{N}
  +\frac{\log\left|J_N(E,\omega)\right|}{N},
\end{align*}
On the other hand
\begin{align*}
  \lim_{N\to\infty}
  \frac{\log\left|J_N(W;\omega)\right|}{N}
  &\underset{?}{=}
  \lim_{r\to1}
  \left\{
    \lim_{N\to\infty}
    \frac{\log\left|J_N(W;\omega^r)\right|}{N}
  \right\}
  \\
  &\underset{?}{=}
  \lim_{r\to1}
  \left\{
    \frac{\log\left|v(N;\omega^r)J_{N,N}(H;\omega^r)\right|}{N}
  \right\}
  \\
  &=
  \lim_{r\to1}\frac{\log\left|v(N;\omega^r)\right|}{N}
\end{align*}
since $J_{N,N}(W,\omega^r)$ can be expressed in terms of sine of $1/N$.
Therefore if we accept these calculations, we could prove
\begin{equation*}
  \lim_{n\to\infty}\frac{\log\left|J_N(D(E),\omega)\right|}{N}
  =
  \lim_{n\to\infty}\frac{\log\left|J_{N,N}(W,\omega)\right|}{N}
  +
  \lim_{n\to\infty}\frac{\log\left|J_{N}(E,\omega)\right|}{N}.
\end{equation*}
Noting that the complement of $D(E)$ is the union of those of the figure-eight
knot and the Whitehead link, which is the volume conjecture for the
Whitehead double of the figure-eight knot.
\begin{ack}
This article is prepared for the proceedings of the workshop
`Volume Conjecture and Its Related Topics' held at the International Institute
for Advanced Studies from 5th to 8th March, 2002.
It was also supported by the Research Institute for Mathematical
Sciences, Kyoto University.
I would like to thank the institutes for their hospitality.
I also thank Kazuhiro Hikami for introducing the computer program PARI.
\par
Part of this work was done when I was visiting Warwick University to
attend the workshop `Quantum Topology' from 18th to 22nd March, 2002,
and Universit{\'e} du Qu{\'e}bec {\`a} Montr{\'e}al to attend the workshop
`Knots in Montreal II' from 20th to 21st April, 2002.
Thanks are due to the universities and to the organizers of the workshops,
Stavros Garoufalidis, Colin Rourke, Steven Boyer, and Adam Sikora.
\end{ack}
\bibliography{mrabbrev,hitoshi}
\bibliographystyle{amsplain}
\end{document}